# Extending a Limit-Free Algebraic–Geometric Construction of the Derivative to Elementary Functions

**Davit Kapanadze**
Doctor of Pedagogical Sciences, Professor
Georgian National University (SEU), Faculty of Business and Technology
Associate Member of the Institute of Mathematics and its Applications (IMA, UK)
ORCID: 0009-0006-8112-5084

## Abstract

The present paper is a continuation of the author's previously published study, in which an algebraic–geometric construction of the derivative was developed in the class of polynomial functions without the initial use of limits:
Kapanadze, D. (2026).
*A limit-free algebraic–geometric construction of the derivative in the class of polynomial functions.*
arXiv:2604.20888.
In the present article, the proposed approach is extended to a broader class of elementary functions. The derivatives of rational power, exponential, logarithmic, and trigonometric functions are considered and obtained on the basis of the geometric interpretation of the tangent, the symmetry of inverse functions, and the local linear structure of a function.
The central idea of the paper is to regard the derivative from the outset as a functional correspondence assigning to each point the slope coefficient of the tangent line. The proposed approach aims to reduce formal abstraction at the introductory stage of differential calculus and to strengthen the conceptual understanding of the derivative. Within this framework, geometric intuition precedes analytic formalisation: the tangent, motion, and local direction of a function are considered first, while the limiting apparatus appears as an analytic continuation of an already constructed geometric–algebraic model.
Thus, the paper does not reject classical limit-based analysis; rather, it proposes a different methodological sequence:

**Tangent → Local Linear Structure → Limit Formalisation.**

From this perspective, the proposed approach may be regarded as a conceptual bridge between geometric intuition, algebraic construction, and classical analysis.

**Keywords**
Derivative without limits; tangent line; algebraic–geometric method; elementary functions; local linear structure; differential calculus; methodology of mathematics education; inverse functions; uniqueness of the tangent.

## Introduction

The concept of the derivative constitutes one of the fundamental foundations of mathematical analysis. It describes the rate of change of a function, the direction of motion, and the slope of the tangent line. In classical analysis, the derivative is defined through the limit of the ratio of increments, which ensures the formal rigour and logical completeness of the theory. However, at the introductory stage of instruction, this very formalism often becomes one of the major conceptual difficulties for pupils and students.

Within the traditional approach, the learner is confronted simultaneously with several abstract ideas: infinitesimal change, limiting transition, formal symbolic notation, and the geometric interpretation of the derivative. As a result, the derivative is frequently perceived not as a natural description of functional change, but rather as the outcome of mechanical manipulations with formulas. At the same time, practical problems involving motion, optimisation, and the determination of maxima and minima are generally perceived far more intuitively.

These circumstances motivated the development of the present algebraic–geometric approach, whose principal idea is to construct the derivative on the basis of the existence and uniqueness of the tangent line, without treating the limit as the initial logical axiom. In the proposed framework, the derivative is introduced from the outset as a function assigning to each point the corresponding slope coefficient of the tangent line.

In the author's previously published work devoted to the algebraic–geometric construction of the derivative in the class of polynomial functions, a criterion for the existence and uniqueness of the tangent line was established. According to this criterion, if the difference between a function and a corresponding linear function is divisible by $(x-a)^2$ at a given point, then the corresponding line is a tangent line. This model demonstrated that the derivative may be constructed from the geometric idea of tangency without the direct initial use of limit formalisation.

The present paper extends this approach to a broader class of elementary functions. Rational power, exponential, logarithmic, and trigonometric functions are considered, and their derivatives are obtained through inverse symmetry, local linear structure, and the analysis of geometric motion.

The central idea of the proposed approach is the local structure of a function. If a function $f(x)$ possesses a tangent line at the point $x_0$, then its local representation may be written in the form

$$f(x) = f(x_0) + k(x - x_0) + R(x),$$

where $k$ denotes the slope coefficient of the tangent line, while $R(x)$ represents a higher-order infinitesimal remainder.

The chapters of the present paper develop this algebraic–geometric approach progressively across different classes of elementary functions. In the case of rational power functions, the symmetry of inverse functions is employed; for exponential and logarithmic functions, the idea of local growth is used; and for trigonometric functions, the geometry of the unit circle and the direction of motion provide the fundamental framework.

Thus, the paper demonstrates that the algebraic–geometric idea of tangency constructed for polynomial functions extends naturally to a wider class of elementary functions. From this standpoint, the derivative may be regarded not as an initially abstract definition, but as a geometric–algebraic model of functional change whose analytic formalisation ultimately becomes consistent with classical limit theory.

## Chapter I
## Derivative of a Rational Power Function
## (Symmetry of the Tangent and Local Structure)

### 1.1 Formulation of the Problem

Following the algebraic–geometric construction of the derivative in the class of polynomial functions, a natural question arises: can a similar idea be extended to elementary functions?
In the present chapter, we consider the rational power function

$$f(x) = x^{p/q},$$

where $p, q \in \mathbb{N}$.
The principal objective is to construct the derivative on the basis of the geometric idea of the tangent line, the symmetry of inverse functions, and the local structure of a function, in such a way that the limiting process does not serve as the initial constructive mechanism.
From a methodological perspective, this approach is significant because the derivative is perceived from the outset as a functional description of the slope of the tangent line rather than merely as the result of a formal limiting operation. This contributes to a clearer understanding of geometric intuition and the local behaviour of functions.

### 1.2 The Initial Case — The Root Function

Consider the function

$$y = \sqrt[q]{x} = x^{1/q}.$$

This function is the inverse of the polynomial function

$$x = y^q.$$

In the previous work, the following result was obtained:

$$(y^q)' = qy^{q-1}.$$

Therefore, the derivative of the root function may be analysed by using the idea of tangent symmetry.
Geometrically, this means that if two functions are inverse to one another, then their graphs are symmetric with respect to the line

$$y = x,$$

while the slope coefficients of the corresponding tangent lines are reciprocal quantities.
Such an approach allows the learner to recognise that a new result often emerges as a natural continuation of an already known structure rather than as a completely independent formula.

### 1.3 Principle of Inverse Functions
### (Symmetry of Tangents)

If the functions

$$y = f(x), x = g(y)$$

are inverse to one another, then their graphs are symmetric with respect to the line

$$y = x.$$

As a consequence of this symmetry, the slope coefficients of the tangent lines at corresponding points are reciprocal:

$$g'(y) = \frac{1}{f'(x)}.$$

Geometrically, this fact means that the tangent line to the inverse function is obtained through the mirror transformation of the original tangent line.
At this stage, the derivative is interpreted not merely as an algebraic operation, but as a local structure connected with geometric symmetry.

### 1.4 Construction of the Derivative of the Root Function

Suppose that the function

$$x = y^q$$

is given.
Then

$$\frac{dx}{dy} = qy^{q-1}.$$

Since the inverse function is

$$y = x^{1/q},$$

the principle of the derivative of the inverse function yields

$$\frac{dy}{dx} = \frac{1}{qy^{q-1}}.$$

Substituting

$$y = x^{1/q},$$

we obtain

$$\frac{dy}{dx} = \frac{1}{qx^{(q-1)/q}} = \frac{1}{q}x^{1/q-1}.$$

Therefore,

$$\left(x^{1/q}\right)' = \frac{1}{q}x^{1/q-1}.$$

The obtained formula is naturally connected with the already known polynomial structure and with the idea of tangent symmetry. As a result, the differentiation rule is perceived as a continuation of logically interconnected processes.

### 1.5 The General Rational Power

Now consider the general case

$$f(x) = x^{p/q}.$$

Observe that

$$x^{p/q} = (x^{1/q})^p.$$

Let

$$u = x^{1/q}.$$

Then the function takes the form

$$f(x) = u^p.$$

Thus, a rational power function may be regarded as a combination of two interconnected structures:

- the root function;
- and the polynomial power.

This approach makes it possible to represent a more complicated function as a natural combination of already known structures.

### 1.6 Principle of the Composite Function

Since

$$f(x) = u^p,$$

we have

$$(u^p)' = pu^{p-1}.$$

And since

$$u = x^{1/q},$$

it follows that

$$u' = \frac{1}{q}x^{1/q-1}.$$

Therefore,

$$f'(x) = pu^{p-1} \cdot u'.$$

Substituting

$$u = x^{1/q},$$

we obtain

$$f'(x) = \frac{p}{q} x^{p/q-1}.$$

Hence,

$$\left(x^{p/q}\right)' = \frac{p}{q} x^{p/q-1}.$$

The present construction is based on the interrelation between root and polynomial structures together with the local behaviour of the function. Consequently, the rules of differentiation are perceived as interconnected processes rather than as independently memorised formulas.

### 1.7 Final Result and Its Interpretation

The classical power rule has been obtained:

$$\left(x^{p/q}\right)' = \frac{p}{q} x^{p/q-1}.$$

Within the proposed approach, this formula is constructed on the basis of:

- the idea of the existence and uniqueness of the tangent line;
- the symmetry of inverse functions;
- and the local structure of the function.

From this perspective, the derivative is regarded as a geometric description of the local change of a function.

In the learning process, such an approach enables the learner to connect the derivative with motion, direction, and the intuitive picture of change, after which analytic formalisation acquires a more natural character.

### 1.8 Special Case — The Function $f(x) = |x|$

The function

$$f(x) = |x|$$

represents an important example because it clearly demonstrates the connection between the non-existence of a tangent line and non-differentiability (Figure 1.1).

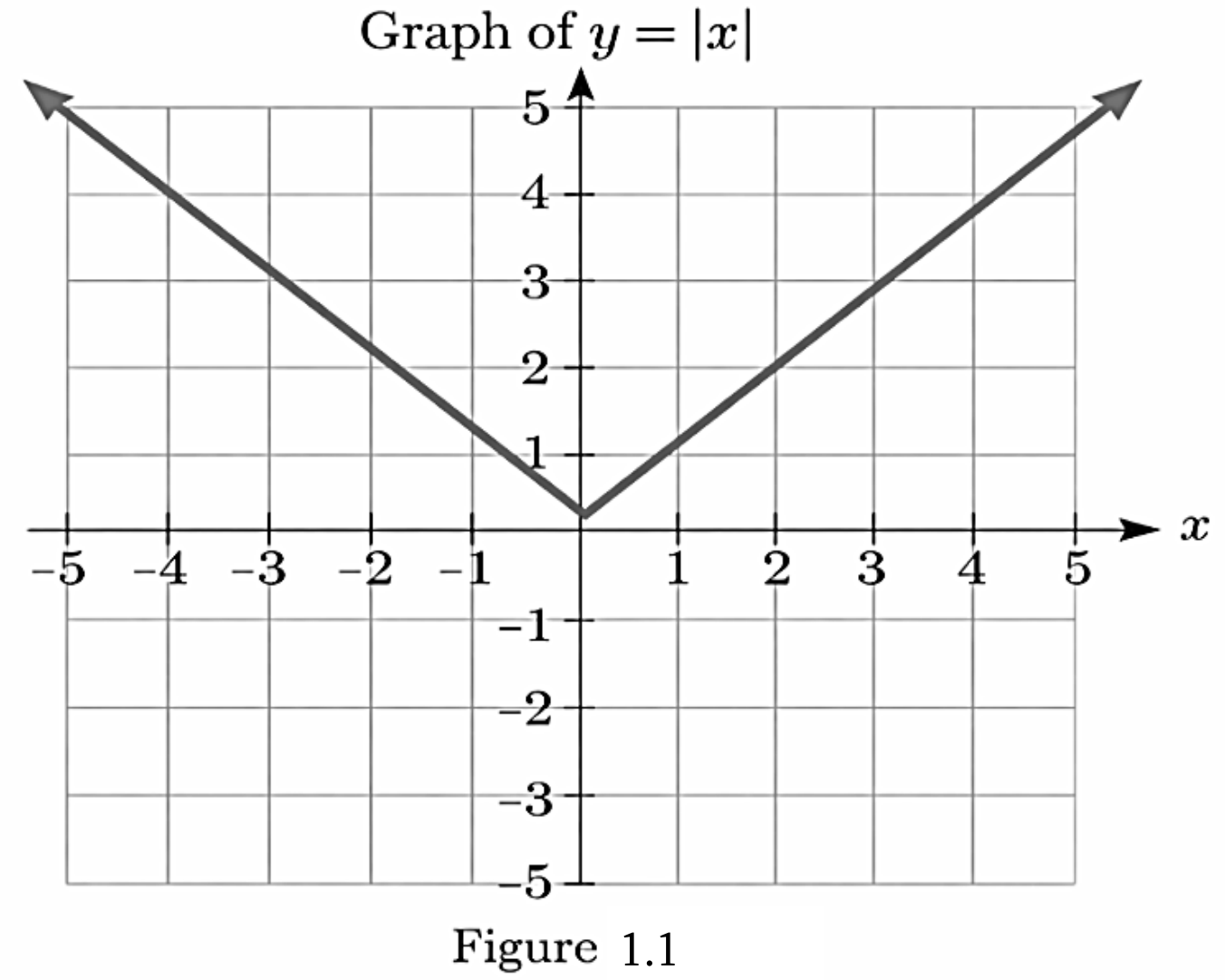

Figure 1.1

At the point $x = 0$, the graph of the function has a corner point.
From the left-hand side, the slope coefficient of the tangent line is

$$k_{-} = -1,$$

whereas from the right-hand side,

$$k_{+} = 1.$$

Thus, at the point

$$x = 0,$$

there exist two different directions, and the tangent line is no longer unique.

### 1.8.1 Algebraic Diagnosis of Non-Differentiability

According to the proposed algebraic criterion, the derivative exists when it is possible to find such a value of $k$that

$$f(x) - (f(x_0) + k(x - x_0))$$

is divisible by

$$(x - x_0)^2.$$

In the present case:

- from the left-hand side, $k = -1$;
- from the right-hand side, $k = 1$.

Since no single unified value exists, the uniqueness of the tangent line is violated, and the algebraic criterion cannot be satisfied.

In this example, non-differentiability is perceived not merely as a formal fact, but as a breakdown of the local structure of the function.

### 1.8.2 Conclusion

At the point $x = 0$,

- the tangent line is not unique;
- the algebraic criterion is not satisfied;
- therefore, the derivative does not exist.

Thus, the proposed approach makes it possible to regard the idea of differentiability as a question concerning the existence and uniqueness of the local tangent structure.

# Chapter II
## Derivative of the Exponential Function
## (Local Growth and the Geometry of the Tangent Line)

### 2.1 Formulation of the Problem

In the previous chapter, the derivatives of rational power functions were considered by means of the symmetry of inverse functions and local tangent structures. We now proceed to exponential functions, whose behaviour plays a particularly important role both in mathematical analysis and in the modelling of natural processes.
Consider the function

$$y = e^x.$$

This exponential function represents one of the fundamental mathematical models of growth. Population growth, radioactive decay, financial interest, and many physical and biological processes are described precisely through exponential dependence.
In classical analysis, the derivative of the function is obtained through a limiting process. In the present approach, the principal objective is to describe the slope coefficient of the tangent line on the basis of the local behaviour and geometric structure of the function, without treating the limit as the initial constructive mechanism.

### 2.2 Local Structure of the Exponential Function

Consider the function $y = e^x$. (Figure 2.1).

The function is positive and increasing for all values of $x$. Its essential property is that the structure of local change preserves the same form at every point.
Let us take the point

$$A(x_0, e^{x_0})$$

on the graph.

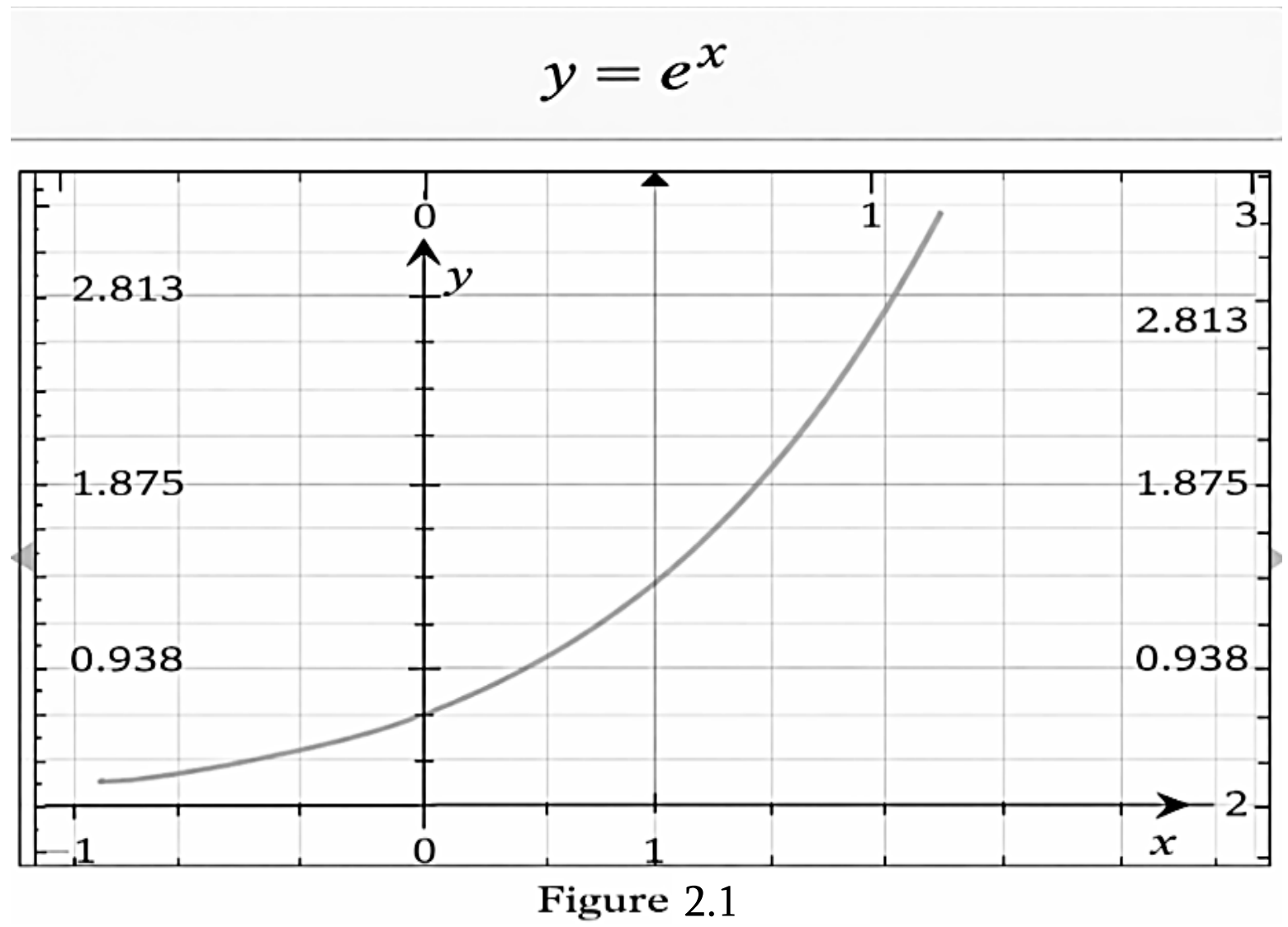

Figure 2.1

Assume that the slope coefficient of the tangent line at this point is $k$.
Then the equation of the tangent line is

$$y = e^{x_0} + k(x - x_0).$$

Thus, the problem is reduced to finding such a value of $k$ for which the line locally touches the graph of the function and establishes first-order contact with it.
Within this approach, the formula appears as an analytic description of the local geometric behaviour of the function.

### 2.3 Geometric Justification of the Tangent Line

Let us use the well-known elementary inequality

$$e^t \geq 1 + t,$$

where equality holds only when

$$t = 0.$$

Geometrically, this fact means that the corresponding line touches the graph of the function at only one point and does not locally intersect it.
Let us write

$$t = x - x_0.$$

Then we obtain

$$e^{x - x_0} \geq 1 + (x - x_0).$$

Multiplying both sides by $e^{x_0}$, we get

$$e^x \geq e^{x_0} + e^{x_0}(x - x_0).$$

From this it follows that the line

$$y = e^{x_0} + e^{x_0}(x - x_0)$$

lies on one side of the graph and touches it only at the point $x_0$. (Figure 2.2).

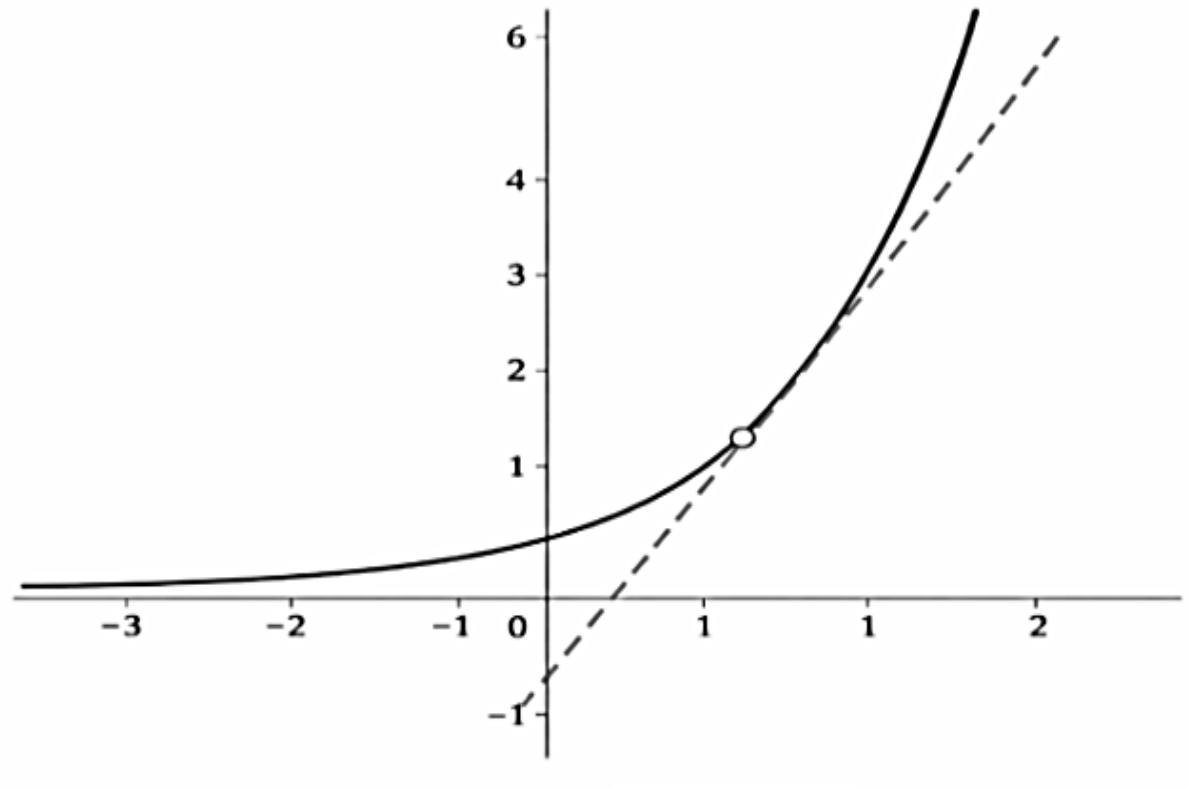

**Figure** 2.2

Therefore,

$$k = e^{x_0}.$$

This result shows that the rate of local change of the exponential function coincides with the value of the function itself.

**2.4 Final Result**

Hence, the slope coefficient of the tangent line is

$$e^{x_0}.$$

Since $x_0$ is an arbitrary point, we obtain

$$(e^x)' = e^x.$$

Thus, the derivative of the exponential function is described on the basis of the geometric criterion of tangency and the local structure of the function.
At every point, the exponential function preserves its own pattern of growth, which determines its special role in mathematics and the applied sciences.

### 2.5 Methodological Interpretation

The proposed approach connects the derivative with the local direction, growth, and rate of change of a function.
From this perspective:

- the tangent line is interpreted as a geometric model of the local behaviour of a function;
- while the derivative is regarded as a description of the intensity of change.

Such a sequence enables the learner to perceive differential calculus not merely as a formal symbolic technique, but as an instrument for the mathematical description of real processes.

## Chapter III
## Exponential and Logarithmic Functions
## (Inverse Symmetry and the Functional Nature of the Tangent Line)

### 3.1 Preliminary Considerations

In the previous chapters, the derivatives of rational power and exponential functions were obtained on the basis of the geometric idea of the tangent line and local linear behaviour. The present chapter extends this approach to general exponential and logarithmic functions.
We consider:

- the function $a^x$, where $a > 0$;
- the natural logarithm;
- and logarithms with an arbitrary base.

The fundamental idea of the proposed method remains unchanged:
the derivative is regarded as the function describing the slope of the tangent line, while the relationship between functions is interpreted as a manifestation of local geometric symmetry.
Logarithmic and exponential functions are especially significant examples because they are inverse to one another. For this reason, the tangent symmetry between them reveals the conceptual unity of the proposed method with particular clarity.
From a methodological perspective, this stage is also important because the learner begins to perceive deep connections between functions. Formulas are no longer viewed as isolated rules; rather, they appear as parts of a unified geometric–algebraic structure. This strengthens the perception of mathematical coherence and reduces the need for purely formal memorisation.

### 3.2 Geometry of Inverse Symmetry

Two functions $f$ and $g$ are inverse to one another if

$$g(f(x)) = x.$$

Their graphs are symmetric with respect to the line $y = x$. (Figure 3.1).

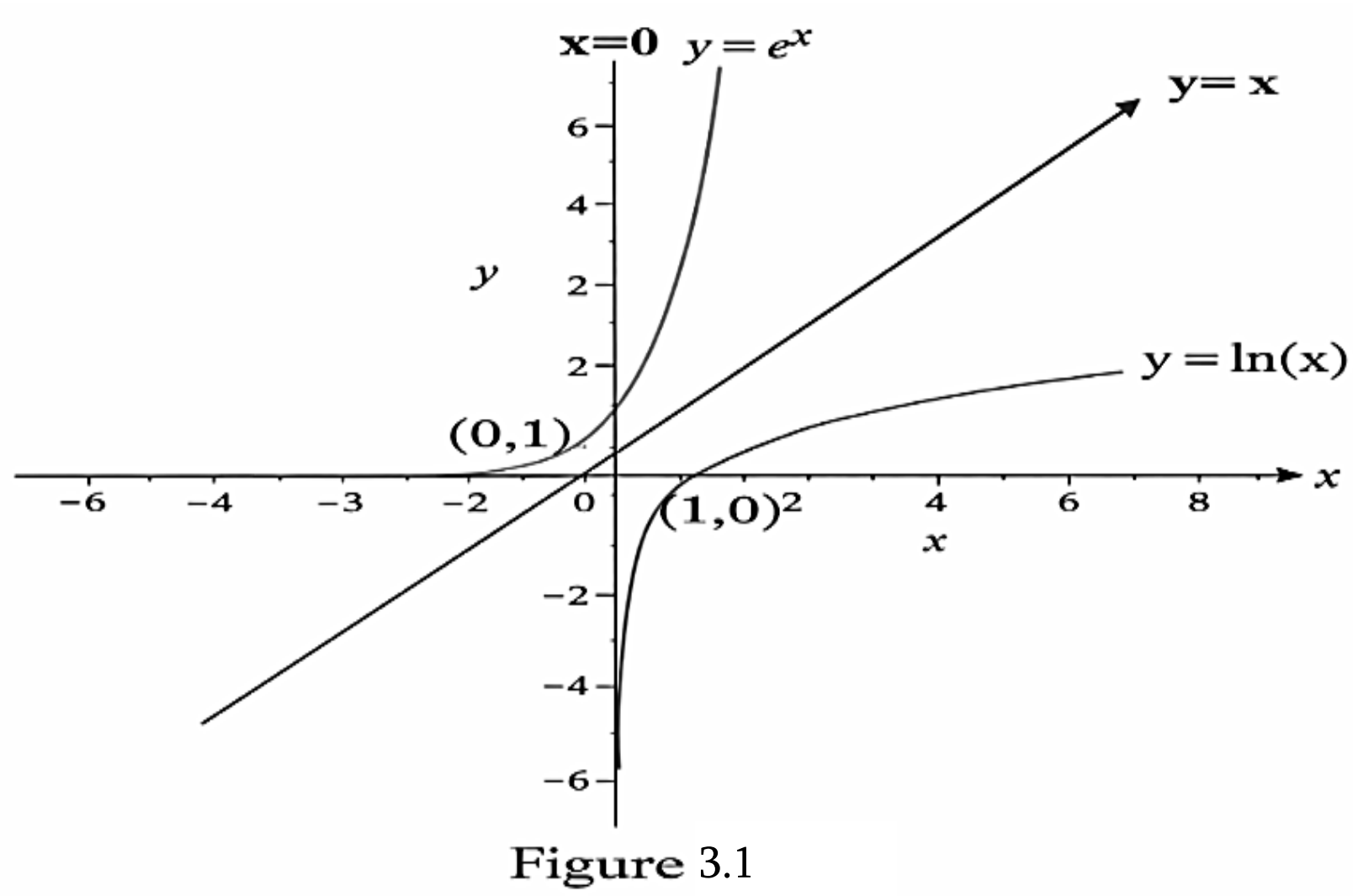

Figure 3.1

If the point $(x_0, y)$ belongs to the graph of $f$, then the point $(y, x_0)$ belongs to the graph of $g$.
As a consequence of this geometric symmetry, the slope coefficients of the tangent lines at corresponding points are reciprocal.
Thus,

$$g'(y) = \frac{1}{f'(x)}.$$

This principle acquires particular significance in the learning process. The learner begins to recognise that the derivative is not merely a local computation; it is connected with the global geometric structure, motion, and symmetry of a function.
Such intuitive representations substantially reduce the sense of abstraction associated with analysis and contribute to a deeper conceptual understanding of mathematical ideas.

### 3.3 Derivative of the Inverse Function

In order to apply the derivative of the inverse function, it is necessary that:

- the function be single-valued;
- the function be monotonic;
- and $f'(x) \neq 0$.

Under these conditions, the following relation holds:

$$g'(y) = \frac{1}{f'(x)},$$

where

$$y = f(x).$$

Geometrically, this means that if the tangent line to the original function moves with a certain slope, then the tangent line to the inverse function acquires the mirror-reciprocal form of this slope.
Thus, the derivative of the inverse function is obtained not as an independent formula, but as a natural consequence of geometric symmetry. This approach emphasises the structural connection between functions and integrates the derivative into a unified geometric model.

### 3.4 Derivative of the General Exponential Function

Consider the function

$$y = a^x, a > 0.$$

Let us use the well-known representation

$$a^x = e^{x \ln a}.$$

Since the derivative of the exponential function $e^x$

has already been established, the derivative of the new function may be constructed through a composite structure.
We have:

$$\frac{d}{dx} e^{x \ln a} = e^{x \ln a} \cdot \ln a.$$

Therefore,

$$\frac{d}{dx} a^x = a^x \ln a.$$

This result clearly shows that the rate of change of an exponential function depends both on the function itself and on the base of growth $a$.
Thus, the derivative of the general exponential function emerges as a natural continuation of already established local structures rather than as an isolated rule.

### 3.5 Derivative of the Natural Logarithm

The function

$$y = e^x$$

is strictly increasing on the whole of $\mathbb{R}$, and therefore it is invertible.
Its inverse function is

$$x = \ln y.$$

Since

$$\frac{d}{dx} e^x = e^x,$$

the principle of the inverse function yields

$$\frac{d}{dx}\ln x = \frac{1}{x}.$$

Hence,

$$\left(\ln x\right)' = \frac{1}{x}.$$

The derivative of the logarithmic function is directly connected with the geometric structure of the exponential function. This once again demonstrates that the proposed approach is based on interconnected functional and geometric ideas.

### 3.6 Logarithm with an Arbitrary Base

Since

$$\log_a x = \frac{\ln x}{\ln a},$$

where $\ln a$ is constant, we obtain

$$\frac{d}{dx}\log_a x = \frac{1}{x\ln a}.$$

Therefore,

$$\left(\log_a x\right)' = \frac{1}{x\ln a}.$$

This formula also follows naturally from the already established logarithmic and exponential structures, once again demonstrating the internal coherence of the proposed system.

### 3.7 Methodological Summary

The discussion of exponential and logarithmic functions has shown that their derivatives may be obtained on the basis of:

- the geometric idea of the tangent line;
- inverse symmetry;
- and local linear structure.

Within this approach:

- the derivative acquires a functional nature from the outset;
- geometric intuition precedes formal symbolism;
- while the limit appears as the analytic formalisation of an already constructed structure.

Thus, the derivatives of exponential and logarithmic functions are regarded not as isolated rules, but as natural continuations of a unified geometric–algebraic system.

## Chapter IV
## Derivatives of Trigonometric Functions
## (The Unit Circle and the Geometric Idea of the Tangent Line)

### 4.1 Introduction
In classical mathematical analysis, the derivatives of trigonometric functions are obtained through special limiting processes. In the present approach, the same formulas are derived from the geometric interpretation of the tangent line and the local direction of motion.
The fundamental principle remains unchanged:
the derivative is the unique slope coefficient of the tangent line.
The natural geometric setting for trigonometric functions is the unit circle. It is precisely here that the connection between

- motion;
- direction;
- rotation;
- and the local change of a function

becomes apparent.
From a methodological perspective, this chapter is particularly significant because the learner encounters, for the first time, the idea that the derivative may be connected not only with algebraic formulas, but also with spatial motion and geometric intuition. Under these conditions, trigonometric functions cease to appear as a system of complicated formulas and instead become a natural mathematical description of motion.

### 4.2 The Unit Circle and Parametrisation
Consider the unit circle with equation

$$x^2 + y^2 = 1.$$

If the angle $t$is measured in radians, then the corresponding point on the circle has coordinates

$$(\cos t, \sin t).$$

Here:

- $x = \cos t$is the abscissa;
- $y = \sin t$is the ordinate.

The radian measure possesses particular importance because, in this case, a small angular change coincides with the change in arc length on the unit circle.
This circumstance enables the learner to perceive trigonometric functions as a living geometric process of motion rather than merely as a collection of formulas. Such an intuitive interpretation reduces the level of abstraction and strengthens conceptual understanding.

**4.3 Tangent Line on the Unit Circle**

Consider the point

$$(\cos t, \sin t)$$

on the unit circle.

Its radius vector is

$$(\cos t, \sin t).$$

The tangent line at the given point is perpendicular to the radius. Therefore, the direction vector of the tangent line is obtained by rotating the radius vector through 90°. (Figure 4.1).

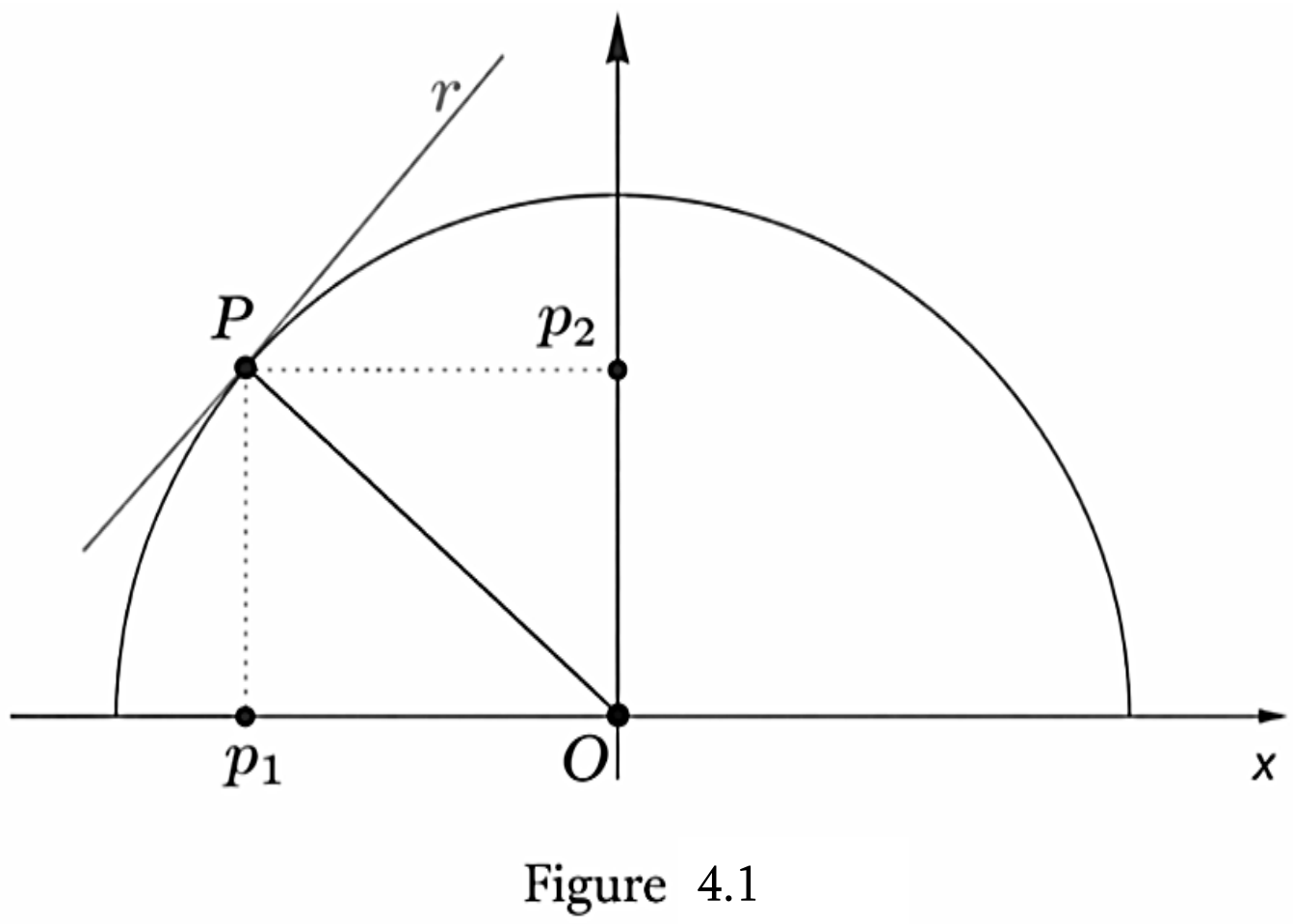

Figure 4.1

Accordingly, we obtain the vector

$$(-\sin t, \cos t).$$

Indeed,

$$(\cos t, \sin t) \cdot (-\sin t, \cos t) = 0,$$

which confirms their orthogonality.

Thus, the direction vector of the tangent line is

$$(-\sin t, \cos t).$$

This vector describes the instantaneous direction of motion along the circle.

In the learning process, this interpretation is particularly effective because the learner directly observes the connection between:

- motion;
- change of direction;
- the tangent line;
- and the derivative.

Such a geometric viewpoint substantially reduces the fear associated with formal symbolism and contributes to an intuitive understanding of the mathematical idea.

### 4.4 Derivative of the Function $\sin x$

The function

$$y = \sin x$$

represents the vertical coordinate of a point on the unit circle.
The vertical component of the tangent direction vector is $\cos x$.

Therefore,

$$\frac{d}{dx}\sin x = \cos x.$$

Here, the derivative naturally appears as the rate of vertical change of motion.
For the learner, this approach is particularly intuitive because the derivative is first perceived as a component of the direction of motion rather than as a formal limiting operation.

### 4.5 Derivative of the Function $\cos x$

The function

$$y = \cos x$$

represents the horizontal coordinate of a point on the unit circle.
The horizontal component of the tangent direction vector is $-\sin x$.

Therefore,

$$\frac{d}{dx}\cos x = -\sin x.$$

At this stage, it becomes especially clear that the variation of trigonometric functions is connected with rotational motion. Consequently, the derivative is perceived as a natural description of the change in direction of motion.

### 4.6 Derivative of the Function $\tan x$

Since

$$\tan x = \frac{\sin x}{\cos x},$$

let us apply the quotient rule:

$$\frac{d}{dx}\tan x = \frac{\cos^2 x + \sin^2 x}{\cos^2 x}.$$

Because

$$\sin^2 x + \cos^2 x = 1,$$

we obtain

$$\frac{d}{dx}\tan x = \frac{1}{\cos^2 x}.$$

At this stage, the learner already sees that the previously established trigonometric and algebraic structures unite into a single coherent system. This strengthens the perception of mathematical unity and supports the development of analytic thinking.

### 4.7 Derivative of the Function $\cot x$

Similarly, we obtain

$$\frac{d}{dx}\cot x = -\frac{1}{\sin^2 x}.$$

Here as well, the derivative is naturally connected with the geometric structure of trigonometric motion.

### 4.8 Methodological and Conceptual Summary

Using the geometric model of the unit circle, the following results were obtained:

$$\frac{d}{dx}\sin x = \cos x,$$
$$\frac{d}{dx}\cos x = -\sin x,$$
$$\frac{d}{dx}\tan x = \frac{1}{\cos^2 x},$$
$$\frac{d}{dx}\cot x = -\frac{1}{\sin^2 x}.$$

The obtained results fully coincide with the formulas of classical analysis; however, within the proposed approach they are derived from:

- the geometric interpretation of motion;
- the direction of the tangent line;
- and the analysis of local behaviour.

From a methodological perspective, this approach provides a significant advantage. The learner:

- observes motion;
- perceives direction;
- grasps the geometric meaning of the tangent line;
- and only afterwards proceeds to formal analytic notation.

Such a sequence reduces the sense of alienation and difficulty often associated with trigonometric analysis, strengthens conceptual understanding, and creates a foundation for perceiving the derivative as a natural mathematical model of real motion and change.

### 4.9 Inverse Trigonometric Functions

Since the principle of tangent symmetry for inverse functions has already been discussed above, the same method may be used to obtain the derivatives of inverse trigonometric functions.
For example, if

$$y = \arcsin x,$$

then

$$x = \sin y.$$

Since

$$\frac{d}{dy}\sin y = \cos y,$$

the principle of the inverse function gives

$$\frac{d}{dx}\arcsin x = \frac{1}{\cos y}.$$

From the unit circle,

$$\cos y = \sqrt{1 - x^2}.$$

Therefore,

$$\frac{d}{dx}\arcsin x = \frac{1}{\sqrt{1 - x^2}}.$$

Similarly, we obtain:

$$\frac{d}{dx}\arccos x = -\frac{1}{\sqrt{1 - x^2}},$$
$$\frac{d}{dx}\arctan x = \frac{1}{1 + x^2},$$
$$\frac{d}{dx}\mathrm{arccot} x = -\frac{1}{1 + x^2}.$$

Thus, the derivatives of inverse trigonometric functions are also obtained naturally from:

- the geometric interpretation of the tangent line;
- inverse symmetry;
- and the local structure of the unit circle.

From a methodological point of view, this approach once again demonstrates that differentiation formulas may be understood as a system of interconnected geometric–algebraic structures rather than as a collection of isolated rules.

## Chapter V
## Algebraic Analysis and the Conceptual Reconstruction of the Derivative

### 5.1 Main Result of the Study

In the present paper, the concept of the derivative has been constructed not from the initial definition based on limiting processes, but from the algebraic–geometric principle of the existence and uniqueness of the tangent line.
The study has progressively carried out:

- the algebraic construction of the tangent line for polynomial functions;
- the extension of the method to rational power functions;
- the construction of the local structure of exponential and logarithmic functions;
- the derivation of the derivatives of trigonometric functions through the geometry of the unit circle;
- and the application of the symmetry of inverse functions.

Within this process, the derivative was regarded from the outset as the functional correspondence

$$x_0 \mapsto f'(x_0),$$

where $f'(x_0)$ denotes the slope coefficient of the tangent line at the corresponding point.
It is important that, within this approach, the derivative is no longer perceived merely as an isolated numerical result at a particular point. From the very beginning, it acquires a functional, geometric, and dynamic character, which substantially facilitates its conceptual interpretation in the learning process.

### 5.2 Local Linear Decomposition and Analytic Formalisation

The central idea of the proposed system is the local linear behaviour of a function.
If a function $f(x)$ possesses a tangent line at the point $x_0$, then its local representation may be written in the form

$$f(x_0 + \Delta x) = f(x_0) + f'(x_0)\Delta x + R(\Delta x),$$

where $R(\Delta x)$represents a higher-order infinitesimal remainder.
Such a representation unifies:

- the geometric idea of the tangent line;
- the local algebraic structure of the function;
- and analytic formalisation.

Particularly important here is the fact that the learner first perceives the local behaviour of a function through geometric intuition, while analytic formalisation appears as a natural continuation of an already constructed structure.

It is precisely this sequence that reduces the level of abstraction and significantly facilitates the introductory understanding of differential calculus.

### 5.3 The Status of the Limit and Logical Sequence

In classical analysis, the derivative is defined by the formula

$$f'(x_0) = \lim_{\Delta x \to 0} \frac{f(x_0 + \Delta x) - f(x_0)}{\Delta x}.$$

In the present system, the logical sequence is different:

**Tangent → Local Linear Decomposition → Limit Formalisation.**

Within this approach, the limit is neither rejected nor replaced. Rather, it appears as the analytic form of an already constructed geometric–algebraic structure.

This circumstance possesses an important methodological effect. The learner no longer encounters, at the initial stage, the high level of abstraction associated with infinitesimal quantities. First, one observes:

- the tangent line;
- motion;
- the direction of change;
- and the local behaviour of the function.

Only afterwards does one proceed to limit formalisation.

Such a sequence substantially reduces:

- the effect of mathematical anxiety;
- alienation from formal symbolism;
- and the sense of the "unnatural difficulty" of analysis.

### 5.4 Comparison Between the Proposed Approach and Classical Analysis

| Aspect | Classical Analysis | Proposed Approach |
|---|---|---|
| Initial idea | Limit | Uniqueness of the tangent line |
| Fundamental mechanism | Limiting process | Algebraic–geometric criterion |
| Intuitive foundation | Infinitesimal changes | Geometric tangent |
| Path to formalisation | Theory of limits | Local linear decomposition |
| Final formulas | $f'(x)$ | $f'(x)$ |

The proposed system does not represent an alternative to classical analysis. On the contrary, it may be regarded as a conceptual and methodological preliminary stage that prepares the learner for a more natural perception of analytic rigour.

### 5.5 Diagnostic Interpretation of Non-Differentiability

One of the important advantages of the proposed method lies in its diagnostic capability.
For example, for the function

$$f(x) = |x|$$

at the point $x = 0$:

- the tangent line is not unique;
- the left-hand and right-hand directions differ;
- and no single local linear decomposition exists.

Therefore, the derivative does not exist.
In this example, the learner not only states formally that "the derivative does not exist," but also directly observes:

- why the uniqueness of the tangent line breaks down;
- why a local linear structure cannot be obtained;
- and how this is connected with the concept of differentiability.

Such an intuitive diagnostic interpretation strengthens the conceptual understanding of the mathematical notion.

### 5.6 Conceptual and Historical Perspective

The history of differential calculus shows that the idea of the derivative was originally connected precisely with geometric and algebraic representations.
In the works of Newton and Leibniz, central importance was attached to:

- motion;
- change;
- the tangent line;
- and the idea of increments.

Later, Cauchy and Weierstrass established the limit-based analytic formalisation that gave calculus its modern rigour.
The present paper demonstrates that an inverse logical path is also possible:

**Tangent → Local Structure → Limit Formalisation.**

From this perspective:

- geometry creates intuition;
- algebra provides the criterion;
- and analysis provides the formalisation.

The unity of these three components yields a conception of the derivative that is simultaneously:

- conceptually clear;
- methodologically accessible;
- and fully compatible with classical analysis.

### 5.7 Final Conclusion

The proposed construction demonstrates that the concept of the derivative may be developed from:

- the geometric idea of the tangent line;
- algebraic structure;
- and the local behaviour of a function;

and only afterwards acquire the formal form of limit analysis.

Such an approach is especially effective at the introductory stage of instruction because it:

- reduces the level of abstraction;
- increases learner motivation;
- facilitates an intuitive understanding of the derivative;
- and strengthens the perception of the interconnection between mathematical ideas.

As a result, the derivative is perceived not as a formal symbolic procedure, but as a natural mathematical language of change, motion, and optimisation processes.

## References

[1] Kapanadze, D. (2026).
*A limit-free algebraic–geometric construction of the derivative in the class of polynomial functions.*
arXiv preprint arXiv:2604.20888.
Available at: https://arxiv.org/abs/2604.20888**erences**

[2] Cauchy, A.-L. (1821). *Cours d'Analyse de l'École Royale Polytechnique.* Paris.

[3] Weierstrass, K. (1894). *Mathematische Werke.* Berlin.

[4] Apostol, T. M. (1974). *Mathematical Analysis* (2nd ed.). Addison–Wesley.

[5] Rudin, W. (1976). *Principles of Mathematical Analysis* (3rd ed.). McGraw–Hill.

[6] Shisha, O. (1986). "The Derivative without Limits."
*Journal of Mathematical Analysis and Applications,* 118, 365–379.

[7] Sangwin, C. J. (2005). "Derivatives without Limits."
*MSOR Connections,* 5(2), 45–50.

[8] Zhang, X., & Tong, Y. (2018). *A Calculus without Limits.*
arXiv preprint arXiv:1802.03029.
Available at: https://arxiv.org/abs/2604.20888

[9] Tall, D. (2010). "A Sensible Approach to the Calculus." Paper presented at the National and International Meeting on the Teaching of Calculus, Pueblo, Mexico.

[10] Robinson, A. (1966). *Non-standard Analysis.* North-Holland.

[11] Keisler, H. J. (1986). *Elementary Calculus: An Infinitesimal Approach.* Prindle, Weber & Schmidt.

[12] Kapanadze, D. (2000). *Foundations of Mathematical Analysis and the School Mathematics Course (A Historical–Methodological Review).*
*Physics and Mathematics at School,* No. 116.
Available online at Academia.edu.
(in Georgian)